\newtheorem{thm}{Theorem}

\newtheorem{myle}{Lemma}

\newtheorem{crl}{Corollary}
\newtheorem{rmk}{Remark}
\documentclass[a4paper, 11pt,english]{article}
\usepackage{amsfonts, graphicx}
\usepackage[iso-8859-7]{inputenc}

\def\i{{\bf i}}
\def\j{{\bf j}}
\def\k{{\bf k}}

\newfont{\Bdd}{msbm10 scaled\magstep1}
\newfont{\footnotesizeBdd}{msbm8 scaled\magstep1}
\begin{document}
 \title{\bf Complex Roots of Quaternion Polynomials}
\author{Petroula Dospra \  and \   Dimitrios Poulakis}
\date{}
\maketitle

\begin{abstract}  The polynomials with quaternion coefficients
 have two kind of roots: isolated and spherical. A spherical root generates a class of roots which contains only
 one complex number $z$ and its conjugate $\bar{z}$, and  this class can be determined by $z$.
In this paper, we deal with the complex roots of quaternion polynomials.
 More precisely,  using B\'{e}zout matrices,
 we give necessary and sufficient conditions,  for a quaternion polynomial to have a complex root, a spherical root,
and  a complex isolated root. Moreover, we compute a bound for the size of the roots of a quaternion polynomial.

\ \\
{\it Keywords: } Quaternion polynomial; B\'{e}zout Matrices;  Spherical Root; Isolated Root.
\ \\
{\it MCS 2010:} 12E15, 11R52, 16H05.
\end{abstract}

\section{Introduction}

Let ${\mathbb R}$  and ${\mathbb C}$ be the fields of real and complex numbers, respectively.
We denote by ${\mathbb H}$ the skew  field of real quaternions. Its elements are of the form
$q =x_0+x_1 \i+x_2 \j+x_3 \k$,
where $x_1,x_2,x_3,x_4\in {\mathbb R}$, and $\i$, $\j$, $\k$ satisfy the
multiplication rules
$$
\i^2= \j^2= \k^2 = -1, \quad
\i \j=- \j \i= \k, \quad
\j \k=- \k \j= \i, \quad
\k \i=- \i \k= \j.
$$
The {\em conjugate}  of $q$ is  defined as $\bar{q} = x_0-x_1 \i-x_2 \j-x_3 \k$. The {\em  real} and the {\em imaginary part}
 of $q$ are  ${\rm Re} q = x_0$ and  ${\rm Im} q = x_1 \i+x_2 \j+x_3 \k$, respectively.
The  {\em norm}  $|q|$  of $q$  is  defined  to  be the quantity
$$|q| = \sqrt{q \bar{q}} = \sqrt{x_0^2+x_1^2+x_2^2+x_3^2}.$$

Two quaternions $q$ and $q^{\prime}$
are said to be {\em congruent} or {\em equivalent,} written $q \sim  q^{\prime}$,
 if for some quaternion $w \neq 0$ we have $q^{\prime} = wqw^{-1}$.   By \cite{Zhang}, we have
$q \sim  q^{\prime}$ if and only if ${\rm Re} q = {\rm Re} q^{\prime}$ and $\ |q| = |q^{\prime}|$.
The  {\em congruence class} of $q$ is the set
$$[q]=\{q^{\prime} \in {\mathbb H}/\ q^{\prime}\sim  q\} =
\{q^{\prime} \in {\mathbb H}/\ {\rm Re} q = {\rm Re} q^{\prime},\ |q| = |q^{\prime}|\}.$$
Note that every class $[q]$  contains exactly one complex number $z$ and its conjugate $\bar{z}$,
which are $x_0 \pm \i \sqrt{x_1^2+x_2^2+x_3^2}$.

 Let ${\mathbb H}[t]$ be the  polynomial  ring  in  the  variable  $t$  over  ${\mathbb H}$.
 Every  polynomial  $f(t)\in  {\mathbb H}[t]$ is  written  as  $a_0t^n+a_1t^{n-l}+\cdots +a_n$
where $n$ is an integer $\geq 0$ and $a_0,\ldots,a_n \in {\mathbb H}$ with $a_0\neq 0$.
The  addition  and  the  multiplication  of polynomials  are  defined  in  the  same  way  as
the  commutative  case,  where  the  variable  $t$  is  assumed  to  commute  with  quaternion
 coefficients \cite[Chapter 5, Section 16]{Lam}.
For every $q \in {\mathbb H}$ we define the evaluation of $f(t)$ at $q$ to be the element
$$f(q) = a_0q^n+a_1q^{n-l}+\cdots +a_n.$$  Note that the evaluation at $q$ is not in general a ring homomorphism from
${\mathbb H}[t]$ to ${\mathbb H}$.

We say that a quaternion $q$ is a zero or a root of $f(t)$ if $f(q) =0$.
In 1941, Niven   proved that the ``Fundamental Theorem of Algebra" holds for
quaternion polynomials \cite{niven}, and in 1944 Eilenberg and I. Niven proved this result for general
quaternion polynomials \cite{eilenberg}. Recent  proofs  can be found in \cite{gentili3} and \cite[Theorem 5.1]{Kalantari}.
The roots of a quaternion polynomial $f(t)$ and its
expression  as a product of linear factors  have been investigated in several papers
\cite{gentili,gentili2,gordon,Janovska,Kalantari,Pogorui,Serodio,Serodio1,Topuridze,Wedderburn}.

A root $q$ of $f(t)$ is called {\em spherical} if $q\not \in {\mathbb R}$ and for every $r \in [q]$ we have
$f(r) = 0$. Otherwise, it is called  {\em isolated.} If two elements of a class  are zeros
of $f(t)$, then  \cite[Theorem 4]{gordon} implies that all elements of this class are zeros of $f(t)$
and so  spherical roots of $f(t)$. Since every class   contains exactly one complex number $z$ and its conjugate $\bar{z}$,
the pairs of complex numbers $\{z,\bar{z}\}$ which are roots of $f(t)$ determine all the spherical roots of $f(t)$.

In this paper we study the complex roots of quaternion polynomials using Bezout matrices.
 First, we determine the degree
of the highest degree complex right factor of a quaternion polynomial, and we give a necessary and sufficient
condition for a quaternion polynomial to have a complex root.
Next, we give necessary and sufficient conditions for a quaternion polynomial to
have a spherical root, and also to have  a complex isolated root.
Finally, we give bounds for the size
of a root of a quaternion polynomial which are sharper than the bound give in \cite[Theorem 4.2]{Opfer}
in case where the roots of polynomial are quite large.

\section{B\'{e}zout Matrices}

In 1971, Barnett computed  the degree (resp.
coefficients) of the greatest common divisor of several univariate polynomials with coefficients in an integral
domain by means of the rank (resp. linear dependencies of the
columns) of several matrices involving theirs coefficients \cite{Barnett,Barnett1}.
 In this section we recall a formulation of Barnett's results using B\'{e}zout matrices \cite{Diaz} which we shall use
 for the presentation of our results. We could equally  use another formulation of Barnett's results given in
\cite{Diaz} or to use another approach \cite{Fatouros,Kakie,Vardulakis}, but we have chosen the formulation with B\'{e}zout matrices as more simple
and quite efficient in computations.

Let $F$ be a field of characteristic zero and
$P (x)$, $Q(x)$ polynomials in $F[x]$ with $d = \max\{\deg P ,\deg Q\}$.
Consider the polynomial
$$\frac{P(x)Q(y)-P(y)Q(x)}{x-y} =  \sum_{i,j=0}^{d-1} c_{i,j}x^iy^j.$$
The {\em B\'{e}zout matrix} associated to $P(x)$ and $Q(x)$ is:
$${\rm Bez}(P, Q) = \left( \begin{array}{cccccccc}
 c_{0,0}      & \cdots & c_{0,d-1}   \\
 \vdots       & \vdots &  \vdots    \\
 c_{d-1,0}    & \cdots & c_{d-1,d-1}
\end{array}
 \right).
$$
The {\em Bezoutian} associated to $P(x)$ and $Q(x)$ is defined as the determinant of the matrix
${\rm Bez}(P, Q)$ and it will be denoted by ${\rm bez}(P, Q)$.
Let $n = \deg P$, $m= \deg Q$ and $p_0$ the leading coefficient of $P(x)$. If $n \geq m$, then
$${\rm bez}(P, Q)=(-1)^{n(n-1)/2} p_0^{n-m} R(P,Q),$$
where $R(P,Q)$ is the well known Sylvester resultant of $P(x)$ and $Q(x)$ \cite{Barnett1,Helmke}.
 Furthermore, we have ${\rm bez}(P, Q)= 0$ if and only if  $\deg(\gcd(P, Q)) \geq 1$.

Now, let $P (x)$, $Q_1(x), \ldots , Q_k (x)$ be a family of polynomials in $F[x]$ with $n = \deg P$ and
$ \deg Q_j \leq n - 1$ for every $j \in \{1, \ldots , k\}$. Set
$${\cal B}_P(Q_1,\ldots,Q_k) = \left( \begin{array}{c}
 {\rm Bez}(P, Q_1)          \\
 \vdots          \\
 {\rm Bez}(P, Q_k)
\end{array}
 \right).
$$

We have the  following formulation of Barnett's theorem.

\begin{myle} The degree of the greatest common divisor of polynomials $P (x)$, $Q_1(x), \ldots , Q_k (x)$
verifies the following formula:
$$\deg(\gcd(P, Q_1, \ldots, Q_k)) = n - {\rm rank}\,{\cal B}_P(Q_1,\ldots,Q_k).$$
\end{myle}
{\it Proof.}\  See \cite[Theorem 3.2]{Diaz}.

\

Moreover, the matrix ${\cal B}_P(Q_1,\ldots,Q_t)$ can provide the greatest common divisor of
$P (x)$, $Q_1(x), \ldots , Q_k (x)$ \cite[Theorem 3.4]{Diaz}.

\section{Complex Roots}

In this section we give a necessary and sufficient condition for a quaternion polynomial
to have a complex root and in this case we determine the solutions of a quadratic equation.

Let $Q(t)\in {\mathbb H}[t]\setminus{\mathbb C}[t]$ be a monic polynomial with $\deg Q = n \geq 1$.
Then there are $f(t), g(t)  \in {\mathbb C}[t]$ with $f(t)g(t)\neq 0$ such that
$$Q(t)= f(t)+\k g(t).$$
We write
 $$f(t) = f_1(t)+f_2(t) \i \ \ \ {\rm and} \ \ \  g(t) = g_1(t)+g_2(t)\i,$$
 where $f_1(t),f_2(t), g_1(t),g_2(t)\in {\mathbb R}[t]$.
 Since $Q(t)$ is monic of degree $n$, we have $\deg f_1 = \deg f = n$ and
$\deg f_2$, $\deg g_1$, $\deg g_2$ are smaller than $n$.

 Set
$$D(t)= \gcd(f_1(t),f_2(t),g_1(t),g_2(t)) \ \ \ {\rm and} \ \ \  E(t)= \gcd(f(t),g(t)).$$
The polynomial $D(t)$ divides $f_1(t)$ and $f_2(t)$, whence we get that $D(t)$ divides $f(t)$.
Similarly, we deduce that $D(t)$ divides $g(t)$. It follows that $D(t)$ divides $E(t)$.

 The polynomial $ B(t)\in {\mathbb H}[t]$ is called a {\em right factor} of
$Q(t)$ if there exists $C(t)\in {\mathbb H}[t]$ such that
$Q(t)=C(t)B(t).$ Note that $q$ is  a root of $Q(t)$ if  and only if $t -  q$ is a right factor of $Q(t)$, i.e. there exists
$g(t) \in  {\mathbb H}[t]$  such  that  $Q(t)  =  g(t)(t -  q)$ \cite[Proposition 16.2]{Lam}. We shall determine
the monic right factors of $Q(t)$ in ${\mathbb C}[t]$ having the highest degree.

\begin{thm} The only  monic right  factor of $Q(t)$ in ${\mathbb C}[t]$ having the highest degree is $E(t)$ and its degree is
$ n - {\rm rank}{\rm Bez}(f,g)$.
 Furthermore, if  $z\in{\mathbb R}$, then $z$ is a root of $ Q(t)$ if and only if it is a root of $D(t)$.
\end{thm}
{\it Proof.}\
Let $G(t)$ be a right factor of $Q(t)$ in ${\mathbb C}[t]$  having the highest degree.
Then there are $a(t),b(t) \in {\mathbb C}[t]$
 such that $Q(t) = (a(t)+\k b(t)) G(t)$. Then $f(t)=a(t)G(t)$ and $g(t)= b(t) G(t)$.
 It follows that $G(t)$ divides $E(t)$. On the other hand, there are $f_1(t),g_1(t)\in {\mathbb C}[t]$ such that
 $f(t) = f_1(t)E(t)$ and $g(t)=g_1(t)E(t)$. Then $Q(t) = (f_1(t)+\k g_1(t)) E(t)$. Since $G(t)$ divides $E(t)$,
 $G(t)$ and $E(t)$ are monic
 and $G(t)$ is a highest degree right factor of $Q(t)$ in ${\mathbb C}[t]$, we deduce that  $G(t) = E(t)$.
 By Lemma 1, we have $\deg E = n - {\rm rank}{\rm Bez}(f,g)$.

Suppose  that $z \in {\mathbb R}$. Then $ Q(z)=0$ if and only if  $f(z) = g(z) = 0$. Since
 $f_1(z),f_2(z),g_1(z),g_2(z)\in {\mathbb R}$, we obtain that
$f(z) = g(z) = 0$ is equivalent to $f_1(z)= f_2(z)=g_1(z)=g_2(z)=0$, and so to $D(z)=0$.

\begin{crl}
The polynomial ${\cal Q}(t)$ has a complex root if and only if we have
$R(f, g)=0$ or equivalently  ${\rm bez}(f, g)= 0$.
\end{crl}
{\it Proof.} \
By Theorem 1,  $ Q(t)$ has a complex root if and only $\deg E > 0$. Further, we have that $\deg E > 0$
if and only $ R(f,g)=0$ which is equivalent to ${\rm bez}(f, g)= 0$.

\begin{crl}
The polynomial $Q(t)$ has at most $n - {\rm rank}{\rm Bez}(f,g)$ complex roots.
\end{crl}

The case of a quadratic quaternion equation has been studied in \cite{Huang,Jia,niven,Dospra}.
 The next corollary provides
theirs solutions in the special case where one of them is complex.

\begin{crl}
Let $Q(t)= t^2+q_1t+q_0$  be a quadratic
polynomial of ${\mathbb H}[t]\setminus {\mathbb C}[t]$ with no real factor.
 Set $q_1= b_1+\k c_1$ and
$q_0 = b_0+ \k c_0$, where $b_0,b_1,c_0,c_1\in {\mathbb C}$.
Then $Q(t)$ has  one complex root if and only if
$$c_0^2- c_0 b_1 c_1+ b_0 c_1^2 = 0.$$
In this case $c_0c_1 \neq 0$,  and the roots of $ Q(t)$ are
 $$q = -\frac{c_0}{c_1}, \ \ \ \ \ \   \sigma = (q-\bar{p})^{-1}p(q-\bar{p}),$$
where $p= - (b_0c_1/c_0+\k c_1)$.
\end{crl}
{\it Proof.}\
Let $f(t) = t^2+b_1t+b_0$ and $g(t)= c_1t+c_0$. We have
$$R(f,g) = c_0^2- c_0 b_1 c_1+ b_0 c_1^2$$
 and by Corollary 1, $Q(t)$ has a complex root if and only if the above quantity is zero.

Suppose now that $ Q(t)$ has a complex root $q$.
If $c_1 = 0$, then the equality $R(f,g)=0$ implies
$c_0 = 0$ and hence $Q(t) \in {\mathbb C}[t]$ which is a contradiction. Thus $c_1 \neq 0$.
If $c_0 = 0$, then we deduce $b_0 =0$, and so $t$ is a factor of $Q(t)$ which is a contradiction.
Therefore  $c_0 \neq 0$.

By Theorem 1, we have $g(q)=0$ and $f(q)=0$. It follows that $q = -c_0/c_1$ and
$f(t) = (t-b_0/q)(t-q)$. Thus, we have the factorization
$$Q(t)= (t-p)(t-q),$$
where $p= - (b_0c_1/c_0+\k c_1)$. If $p = \bar{q}$, then we have
$b_0c_1/c_0+\k c_1 =  \bar{c}_0/\bar{c}_1$.
It follows that $c_1 = 0$ which is a contradiction. Thus, \cite[Lemma 1]{Serodio1} yields
$$Q(t)= (t-(p-\bar{q})q(p-\bar{q})^{-1})(t-(q-\bar{p})^{-1}p(q-\bar{p})).$$
Hence, the other root of $Q(t)$ is
$\sigma = (q-\bar{p})^{-1}p(q-\bar{p}).$

\section{Spherical and Complex Isolated Roots}

In this section, we give necessary and sufficient conditions for a quaternion polynomial to
have a spherical root and to have a complex isolated root.
First, we consider the spherical roots.

\begin{thm}
Let $z \in {\mathbb C}\setminus {\mathbb R}$. The following are equivalent:\\
(a) The number $z$ is a spherical root of $ Q(t)$. \\
(b) The number $z$ and its conjugate $\bar{z}$ are common roots of $f(t)$ and $g(t)$.\\
(c) The number $z$ is a common root of $f_1(t),f_2(t),g_1(t),g_2(t)$.
\end{thm}
{\it Proof.}\  If $z$ is a spherical root of $Q(t)$, then its conjugate $\bar{z}$ is also a root of $Q(t)$.
Thus Theorem 1 implies that $z$ and  $\bar{z}$ are  common roots of $f(t)$ and $g(t)$. If this holds, then
 the polynomial $(t-z)(t-\bar{z})$
is a factor of $f(t)$ and $g(t)$. It follows that $(t-z)(t-\bar{z})$ is a factor of $f_1(t),f_2(t),g_1(t),g_2(t)$.
Hence $z$ is a common root of $f_1(t),f_2(t),g_1(t),g_2(t)$.
Finally, if $z$ is a common root of $f_1(t),f_2(t),g_1(t),g_2(t)$, then $\bar{z}$ is also a common root of
$f_1(t),f_2(t),g_1(t),g_2(t)$. Hence $z$ and $\bar{z}$ are roots of $Q(t)$ and so,
they define the same spherical root.

\begin{crl}
If ${\cal Q}(t)$  has no real factor, then it has only isolated roots.
\end{crl}
{\it Proof.} \  Suppose that $Q(t)$ has a spherical root $\rho$. Then there is a complex number $z \in [\rho]$.
It follows that $z$ is a spherical root of $Q(t)$ and so,  Theorem 2(b) implies that
$z$ and its conjugate $\bar{z}$ are common roots of $f(t)$ and $g(t)$. Thus, the real polynomial $(t-z)(t-\bar{z})$
 is a common factor of  $f(t)$ and $g(t)$. Therefore, $Q(t)$  has the real factor
$(t-z)(t-\bar{z})$ which is a contradiction. Hence $Q(t)$ has only isolated roots.

\begin{rmk} {\rm  Since a spherical root of $ Q(t)$ has in its class a number
 $z \in {\mathbb C}\setminus {\mathbb R}$,  Theorem 2 yields that we can find the spherical roots of $Q(t)$
 by computing all the common complex roots of $f_1(t),f_2(t),g_1(t),g_2(t)$.
}
\end{rmk}

\begin{thm}
Suppose that the quaternion polynomial $ Q(t)$ has no real root. The following are equivalent: \\
(a) The polynomial $Q(t)$ has a spherical root.\\
(b) $\deg D(t) >0$.\\
(c) $n > {\rm rank}\,{\cal B}_{f_1}(f_2,g_1,g_2).$
\end{thm}
{\it Proof.} Suppose that $Q(t)$ has a spherical root $q$. Let $z$  and $\bar{z}$  be the only  complex
numbers of the class of $q$. Then we have $Q(z) = Q(\bar{z})$, whence we get
$$f(z) = f(\bar{z}) = 0 \ \ \  {\rm and }  \ \ \  g(z) = g(\bar{z}) = 0.$$
It follows that the real polynomial $(t-z)(t-\bar{z})$ divides $f(z)$ and $g(z)$ and hence $D(t)$.
Therefore $\deg D(t) > 0$.

Conversely, suppose that $\deg D(t)  > 0$. Then $D(t)$ has a root $z \in   {\mathbb C}$. If $z \in {\mathbb R}$,
then $z$ is a common root of $f_1(t),f_2(t),g_1(t),g_2(t)$ and hence $z$ is a root of $Q(t)$.
Since $ Q(t)$ has no real root we have a contradiction. Thus $z \not \in  {\mathbb R}$ and so, its conjugate
$\bar{z}$  is also a root of $D(t)$.  It follows that $z$ and $\bar{z}$ are roots of $ Q(t)$. By Theorem 2,
the class of $z$ is a spherical root of ${\cal Q}(t)$.

Finally, by Lemma 1 we have
$$\deg D = n - {\rm rank}\,{\cal B}_{f_1}(f_2,g_1,g_2)$$
and so, the equivalence of (b) and (c) follows.

\begin{rmk} {\rm  In the above theorem, the hypothesis that $ Q(t)$ has no real root, implies that
$D(t)$ has not a real root and so, if $\deg D > 0$, then we have that  $\deg D$ is even.
}
\end{rmk}

\begin{thm}
Suppose that the quaternion polynomial $Q(t)$ has no real root.
The following are equivalent:\\
(a) The polynomial $Q(t)$ has an isolated complex root. \\
(b) $\deg E > \deg D$.\\
(c) ${\rm rank}\,{\rm Bez}(f,g)< {\rm rank}\,{\cal B}_{f_1}(f_2,g_1,g_2).$
\end{thm}
{\it Proof.} \
Let $z$ be an isolated complex root of $ Q(t)$. By Theorem 1, $z$ is a common root of $f(t)$ and $g(t)$.
Since the root $z$ is isolated, $\bar{z}$ is not a common root of these two polynomials. Thus, the real polynomial
$(t-z)(t-\bar{z})$ is not a common factor of $f_1(t),f_2(t),g_1(t),g_2(t)$.
Hence $z$ is not a root of $D(t)$.
Since $D(t)$ divides $E(t)$, we deduce that $\deg E > \deg D$.

Conversely, suppose that $\deg E > \deg D$.
Then $E(t)$ has a complex root $z$ which is not a root of $D(t)$. If $z$ is a  spherical root, then $\bar{z}$ is also
a common root of $f(t)$ and $g(t)$. It follows that $(t-z)(t-\bar{z})$ is a common factor of $f_1(t),f_2(t),g_1(t),g_2(t)$
and hence divides $D(t)$. Therefore $z$ is a root of $D(t)$ which is a contradiction. Thus, $z$ is an isolated root of
$Q(t)$.

By Lemma 1, we have
$$\deg D = n - {\rm rank}\,{\cal B}_{f_1}(f_2,g_1,g_2)$$
and
$$\deg E = n - {\rm rank}{\rm Bez}(f,g).$$
Thus, we have
$\deg E >  \deg D$ if and only if
$${\rm rank}{\rm Bez}(f,g)< {\rm rank}\,{\cal B}_{f_1}(f_2,g_1,g_2).$$

\section{Bounds for the Size of the Roots}

In \cite[Section 4]{Opfer} some bounds for the roots of quaternion polynomials are given.
In this section we compute new bounds which are sharper than the bound give in \cite[Theorem 4.2]{Opfer},
in case where the roots of polynomial are quite large.

Let
$$ Q(t) = a_0t^n+a_1 t^{n-1}+\cdots + a_n.$$
be a quaternion polynomial.
We define the {\em height} of $Q(t)$ to be the quantity
 $$H(Q) = \max\{1,|a_1/a_0|,\ldots,|a_n/a_0|\}.$$
We write $Q(t)= f(t)+\k g(t)$, where  $f(t),g(t) \in  {\mathbb C}[t]$, and
 $$f(t) = f_1(t)+f_2(t) \i, \ \ \  g(t) = g_1(t)+g_2(t)\i,$$
 where $f_1(t),f_2(t), g_1(t),g_2(t)\in {\mathbb R}[t]$. We set
$${ \cal H}_1 = \min \{ H(f), H(g)\} \ \  \  {\rm and} \ \ \  { \cal H}_2 = \min \{ H(f_1),H(f_2),H(g_1),H(g_2)\}.$$

\begin{thm}
Suppose that the polynomial $Q(t)$ is monic and $\rho$ is a root of $ Q(t)$. If  $\rho$ is  a spherical root, then
$$|\rho| < 1+{ \cal H}_2^{1/2}.$$
 If  $\rho$ is an isolated complex root, then
 $$|\rho| < 1+{ \cal H}_1.$$
In the general case, we have
$$|\rho| < 1+H( Q).$$
\end{thm}
{\it Proof.}\  Suppose that $\rho$ is a spherical root of $Q(t)$. Then there is
$z \in {\mathbb C}\setminus {\mathbb R}$ in the class of  $\rho$  which is also a root of
 $ Q(t)$. By Theorem 2, $z$ is a  common complex root of $f_1(t),f_2(t),g_1(t),g_2(t)$. Thus,
\cite[Corollary 3]{Mignotte} implies that $|z| < 1+{ \cal H}_2^{1/2}$. Since $|\rho|= |z|$, we obtain
$|\rho| < 1+{ \cal H}_2^{1/2}$.

Suppose that $\rho$ is an isolated root. If $\rho \in {\mathbb C}$, then Theorem 1 implies that
$\rho$ is a common root of $f(t)$ and $g(t)$. Hence \cite[Corollary 2]{Mignotte} yields $|\rho| < 1+{ \cal H}_1$.

Suppose next that $\rho$ is an isolated non-complex root.
If $|\rho| \leq 1$, then the result is true. Suppose that $|\rho| > 1$.
Since $\rho$ is a root of $ Q(t)$, there is $G(t)\in {\mathbb H}[t]$ such that
$Q(t) = G(t) (t-\rho).$
Write
$$G(t) = t^{n-1}+b_1  t^{n-2}+\cdots +b_{n-1}.$$
Then
$$Q(t) =  G(t) (t-\rho)
            =  t^n+(b_1-\rho)t^{n-1}+(b_2-b_1\rho)t^{n-2}+\cdots+b_{n-1}\rho.$$
It follows that
$$a_1 = b_1-\rho, \ \ a_2 = b_2-b_1\rho, \ \ a_3 = b_3-b_2\rho,
  \ldots, \  \  a_n = b_{n-1}\rho.$$
Let $i$ be the smallest index such that $H(G) = |b_i|$. Then we have
$$H( Q)     \geq  |b_i-b_{i-1}\rho |  \geq  ||b_i|-|b_{i-1}\rho ||
 \geq  |H(G)- |b_{i-1}\rho | |
>  |H(G)(1-|\rho |) |,$$
whence we deduce the result.

\begin{rmk}
{\rm
In case where $a_0 = 1$, \cite[Theorem 4.2]{Opfer} yields that the roots $\rho$ of $Q(t)$ satisfy
$$|\rho| \leq \max\{1,\sum_{i=1}^n |a_i| \}.$$
If $\sum_{i=1}^n |a_i| > 1+H(Q)$, then  Theorem 5 gives a sharper bound.}
\end{rmk}

\begin{crl}
Let $Q(t)\in {\mathbb H}[t]\setminus {\mathbb H}$ be a monic polynomial.
Then $ Q(t)$ has at most a finite number of roots $\bf x$
 of the form ${\bf x}= x_1 +x_2  i +x_3 j+  x_4k$,
where $x_1,x_2,x_3, x_4$ are integers.
\end{crl}

{\small
\ \\
Agricultural University of Athens,\\
Department of Natural Resources Management and Agricultural Engineering, \\
Mathematics Laboratory, \\
75 Iera Odos, Athens 11855, Greece\\
pdospra@aua.gr

\ \\ 
Aristotle University of Thessaloniki,  \\
Department of Mathematics, \\
Thessaloniki 54124, Greece, \\
poulakis@math.auth.gr
}

\end{document}